\documentclass[12pt]{amsart}

\newcommand{\allflat}{\mathcal{A}}
\newcommand{\aut}{\sigma}

\newcommand{\all}{\aleph}
\newcommand{\orb}[2]{\mathcal{O}_{#1}(#2)}
\newcommand{\RS}{\Sigma}
\newcommand{\CH}{\mathcal{H}}
\newcommand{\G}{\Gamma}
\newcommand{\C}{\mathbb{C}}
\newcommand{\R}{\mathbb{R}}
\newcommand{\mm}{\mathbf{m}}
\newcommand{\alg}[1]{\mathfrak{#1}}
\newcommand{\forsome}{\text{ for some }}

\newcommand{\Ad}[2][]{\, Ad^{#1} (#2) \, }
\newcommand{\AD}[2][]{\, AD^{#1} (#2) \, }

\newcommand{\SK}{SK}
\newcommand{\mult}[2]{\ensuremath{\mathcal{M}_{#1}(#2)}}
\newcommand{\moduli}[2][\sigma]{\ensuremath{
\mathcal{M}^{#1}_{\RS}(#2)}}

\DeclareMathOperator{\Hol}{Hol}
\DeclareMathOperator{\re}{Re}
\DeclareMathOperator{\im}{Im}

\DeclareMathOperator{\Dim}{dim}
\DeclareMathOperator{\rank}{rank}
\DeclareMathOperator{\Tr}{Tr}

\theoremstyle{plain}
\newtheorem{theorem}{Theorem}
\newtheorem{lemma}{Lemma}

\theoremstyle{definition}
\newtheorem{definition}{Definition}

\theoremstyle{remark}
\newtheorem*{remark}{Remark}

\begin{document}

\title[The hyperbolic moduli space]{The hyperbolic moduli space of
flat connections and the isomorphism of symplectic multiplicity
spaces}

\author{Anton Yu. Alekseev}
\address{Institute for Theoretical Physics \\ Uppsala University \\
Box 803 \\ \mbox{S-75108} Uppsala \\ Sweden}
\email{alekseev@rhea.teorfys.uu.se, malkin@rhea.teorfys.uu.se}

\author{Anton Z. Malkin}
\address{Permanent address of A.M.  is Steklov Mathematical Institute
\\
Fontanka 27 \\ St.Petersburg 191011 \\ Russia}

\date{March 1996}

\begin{abstract}
Let $G$ be a simple complex Lie group, $\alg{g}$ be its Lie algebra,
$K$ be a maximal compact form of $G$ and $\alg{k}$ be a Lie algebra
of $K$. We denote by $X\rightarrow \overline{X}$ the anti-involution
of $\alg{g}$ which singles out the compact form $\alg{k}$.

Consider the space of flat $\alg{g}$-valued connections  on a Riemann
sphere with three holes which satisfy  the additional condition
$\overline{A(z)}=-A(\overline{z})$.   We call the quotient of this
space over the action of the gauge group
$\overline{g(z)}=g^{-1}(\overline{z})$ a
\emph{hyperbolic} moduli space of flat connections.

We prove that the following three symplectic spaces are isomorphic:

1. The hyperbolic moduli space of flat connections.

 2. The symplectic multiplicity space obtained as  symplectic
quotient of the triple product of co-adjoint orbits of $K$.

3. The Poisson-Lie multiplicity space equal to the Poisson quotient
of the triple product of dressing orbits of $K$.

\end{abstract}


\maketitle


\section{Introduction}
In this paper we investigate symplectic multiplicity spaces for
compact simple Lie groups.

It has been recently observed by Jeffrey \cite{Jeffrey} that  in the
case when all three co-adjoint orbits belong to a small vicinity of
$0\in \alg{k}^*$ symplectic multiplicity spaces  coincide with the
moduli spaces of $\alg{k}$-valued flat connections on the sphere with
three holes. It is also known that this result does not hold true for
sufficiently big co-adjoint orbits.

Though we do not know how to extend the results of \cite{Jeffrey}, we
reformulate the problem for the \emph{hyperbolic} moduli space of
flat connections which is obtained as a symplectic quotient of the
space of $\alg{g}$-valued connections with additional condition
$\overline{A(z)}=-A(\overline{z})$ over the action of the gauge group
$\overline{g(z)}=g^{-1}(\overline{z})$. In this setting we establish
the isomorphism of the symplectic multiplicity spaces to the
hyperbolic moduli spaces for \emph{arbitrary} co-adjoint orbits. This
isomorphism is described explicitly.

The construction of the map between  multiplicity spaces and moduli
spaces is based on the simple observation that the holomorphic
connection of the type
\begin{equation}
A(z)=\sum_i \frac{X_i}{z-z_i} dz
\end{equation}
is flat. Here  $X_i$ belong to  given co-adjoint orbits. This
construction has been suggested by  Hitchin in \cite{Hitchin}.

It is worth mentioning that both our main results (subsection
\ref{MultToModuli} and \ref{ModuliToMult}) are inspired by the work
of Fock and Rosly \cite{Fock-Rosly}. As a nice corollary of their
description of the moduli spaces of flat connections we obtain an
isomorphism of the hyperbolic moduli space and the Poisson-Lie
multiplicity space corresponding to an arbitrary Poisson structure on
the compact group $K$.

The paper is organized as follows. In Section \ref{mults} we define
symplectic and Poisson-Lie multiplicity spaces. Section \ref{moduli}
includes a review of the Goldman and Fock-Rosly descriptions of the
moduli space of flat connections. There we also define the hyperbolic
moduli spaces for the sphere with three holes.

The isomorphisms between the multiplicity spaces and the hyperbolic
moduli space
are described  in Section \ref{isomorphisms}.

\section{Symplectic Multiplicity Spaces}\label{mults}

\subsection{Multiplicity spaces for compact Lie groups}

Let $G$ be a complex, connected and simply connected
Lie group with Lie algebra $\alg{g}$ and let $K$ be it's maximal
compact form. We denote by  $\alg{k}^*$
the dual space to the corresponding Lie algebra $\alg{k}$.
There exists a nondegenerate positively definite scalar product on
$\alg{k}$ (Killing form) which we denote by $<,>$. It induces the
canonical isomorphism $I$ between $\alg{k}^*$ and $\alg{k}$:

\begin{equation}
I: \alg{k} \rightarrow \alg{k}^* .
\end{equation}

The group $K$ acts on $\alg{k}^*$ by the coadjoint action $Ad^*$:

\begin{equation}
\Ad[*]{k} =I \circ \Ad{k} \circ I^{-1}.
\end{equation}

The space $\alg{k}^*$ carries a natural Poisson structure
called Kirillov-Kostant Poisson structure which is linear
and invariant under $Ad^*$-action.

\begin{definition}

Kirillov-Kostant bracket
of two functions $\psi$ and $\psi' $ on $\alg{k}^*$ is given
by the following formula:

\begin{equation} \label{KKPb}
\{ \psi ,\psi' \} (P) = <P,[\nabla \psi , \nabla \psi' ]>.
\end{equation}
Here $P \in \alg{k}^*$ and $[,]$ is the commutator
in $\alg{k}$ (we consider $\nabla \psi$
and $\nabla \psi'$ as elements of the Lie algebra $\alg{k}$).

\end{definition}

According to the general theory of Poisson manifolds
\cite{Arnold,Weinstein} the space $\alg{k}^*$ splits into
the set of symplectic leaves. Symplectic leaves
coincide with orbits of the coadjoint action of $K$ on $\alg{k}^*$.
In this paper we consider only coadjoint orbits of the maximal
dimension.
Such orbits are in one to one correspondence
with elements $H$ of the interior of the positive Weyl chamber
$W_+$ in the Cartan subalgebra $ \alg{h} \subset \alg{k}$:
\begin{equation}
\orb{0}{H} =\{ X \in \alg{k}^* : X=\Ad[*]{k} I(H), k \in K, H \in W_+
\} .
\end{equation}
The orbits of  smaller dimension correspond to the elements
$H$ on the boundary of the Weyl chamber.

The main objects of our interest in this paper are
symplectic multiplicity spaces. In order to define them we
consider a direct product of three coadjoint orbits
$\orb{0}{H_1}\times \orb{0}{H_2} \times \orb{0}{H_3}$. This is
a symplectic space with the Poisson
bracket equal to the sum of Kirillov-Kostant brackets
on $\orb{0}{H_1}$, $\orb{0}{H_2}$ and $\orb{0}{H_3}$.

There is a diagonal action of the group $K$ on the triple product :
\begin{multline} \label{kOO}
\orb{0}{H_1}\times \orb{0}{H_2} \times \orb{0}{H_3}
\ni (X_1, X_2, X_3)
\stackrel{k}{\rightarrow}  \\  \stackrel{k}{\rightarrow}
(\Ad[*]{k} X_1 ,\Ad[*]{k} X_2, \Ad[*]{k} X_3 ) .
\end{multline}
This is a Hamiltonian action with the  moment map:
\begin{gather}
\mu :
\orb{0}{H_1} \times \orb{0}{H_2} \times \orb{0}{H_3}
\rightarrow \alg{k}^*:  \nonumber
\\
\mu(X_1,X_2,X_3 )= X_1+X_2+X_3.
\end{gather}

Symplectic multiplicity space is the symplectic quotient of the space
$\orb{0}{H_1} \times \orb{0}{H_2} \times \orb{0}{H_3}$  with respect
to the diagonal  action \eqref{kOO} of $K$.

\begin{definition}
The symplectic multiplicity space $\mult{0}{H_1,H_2,H_3}$
is the quotient  of the space $\mu^{-1}(0)\subset \orb{0}{H_1} \times
\orb{0}{H_2} \times \orb{0}{H_3}$
over the diagonal action \eqref{kOO} of $K$:

\begin{multline}
\mult{0}{H_1,H_2,H_3} =
\\
=\{ (X_1, X_2, X_3) \in
\orb{0}{H_1}\times \orb{0}{H_2} \times \orb{0}{H_3}:
X_1 + X_2 + X_3 = 0 \} /K.
\end{multline}

\end{definition}
We always assume that the stabilizer of a point in $\mu^{-1}(0)$ is a
discrete subgroup of $K$. Then $0$ is said to be a regular value of
the moment map and the quotient space $\mult{0}{H_1,H_2,H_3}$ is a
smooth manifold. As a symplectic quotient it inherits symplectic
structure from the parent space  $\orb{0}{H_1}\times \orb{0}{H_2}
\times \orb{0}{H_3}$ \cite{Weinstein}.

Now we proceed to the definition of Poisson-Lie multiplicity
spaces.

\subsection{Multiplicity spaces for Poisson-Lie groups}\label{PLmult}
The main idea of the theory of Poisson-Lie groups is to equip a Lie
group with a special kind of Poisson structure.

\begin{definition}
A pair $(K, \mathcal{P})$ of  a group $K$ and a Poisson bracket
$\mathcal{P}$ on $K$  is called a Poisson-Lie group if the
multiplication map $K\times K \rightarrow K$ is a Poisson map.
\end{definition}

According to \cite{Lev-Soi} the classes of isomorphic
Poisson-Lie structures on a compact simple Lie group are parametrized
by pairs $(t, u)$ , where $t$ is a real number and $u$ is a real
valued anti-symmetric (with respect to the Killing form) automorphism
of the Cartan subalgebra $\alg{h}$

\begin{gather}
u: \alg{h} \rightarrow \alg{h} \nonumber
\\
<u(X),Y>=-<X,u(Y)>\ , \ \forall X,Y \in \alg{h}.
\end{gather}

In order to describe this family of Poisson brackets we need some
more notations. Let $e_\alpha$ be the  root generator
of $\alg{g}$ corresponding to a root $\alpha$ and $h_i$ be the
generator of
the Cartan subalgebra $\alg{h}$ corresponding to a simple root
$\alpha_i$.
Define a pair of  classical $r$-matrices (solutions to the
classical Yang-Baxter equation \cite{Semenov}) $r^t_+(u) , r^t_-(u)
\in \alg{g} \otimes \alg{g}$, $\alg{g}$ is considered as a Lie
algebra over $\R$:

\begin{multline}
\label{r+t}
r^t_+(u)= \\ =\frac{1}{2}\sum_{j \in \Delta_s} i(h_j+iu(h_j))
\otimes h^j+\sum_{\alpha \in \Delta_+} (e_{\alpha}
\otimes (ie_{\alpha}+ie_{-\alpha }) + ie_{\alpha} \otimes
(e_{-\alpha}-e_{\alpha})),
\end{multline}
\begin{equation}
r^t_{-}(u)=- P(r^t_+(u)),
\end{equation}
where $\Delta_+$ is the set of positive roots, $\Delta_s$ is the set
of simple roots,
$P$ is the flip operation in the tensor product $\alg{g} \otimes
\alg{g}$:
$P(a\otimes b)= b\otimes a$.
We shall often write $r_+^t$ and $r_-^t$ instead of
$r_+^t(u)$ and $r_-^t(u)$. Let us remark that the $r$-matrices
$r_{\pm}^t(u)$ satisfy the following reality condition:
\begin{equation} \label{rreal}
\overline{r^t_+(u)}=r_-^t(u).
\end{equation}
Here $X\rightarrow \overline{X}$ is the anti-involution in $\alg{g}$
which singles out the compact form $\alg{k}$.

\begin{definition}

The Poisson-Lie bracket on the compact group $K$ corresponding to
the pair of classical $r$-matrices $r^t_{\pm}(u)$ is defined as:

\begin{multline} \label{Sb}
\{ \psi , \psi'  \}_{(t,u)} (k)= \\
= <r_{\pm}^t(u), \nabla^L \psi (k) \otimes \nabla^L \psi' (k) > -
<r_{\pm}^t(u), \nabla^R \psi (k) \otimes \nabla^R \psi' (k) >.
\end{multline}
Here   the bracket $<,>$ is  the canonical pairing between
$\alg{k}^*\otimes \alg{k}^*$ and $\alg{k} \otimes \alg{k}$.
\end{definition}
We always understand $\nabla_L$ and $\nabla_R$ as left and right
derivations on the group $G$ considered as a \emph{real} group with
values in the dual space to the Lie algebra $\alg{g}$. Let us note
that the reality condition \eqref{rreal} implies reality of the
Poisson tensor \eqref{Sb}.

Besides the Poisson-Lie group $K$ we shall need the dual Poisson-Lie
group $K^*$ and the Heisenberg double $D$, which provide the
Poisson-Lie counterparts of the dual space $\alg{k}^*$ and of  the
cotangent bundle $T^*K$.

For the Poisson-Lie structures on compact simple Poisson-Lie groups
which we listed above, the dual group $K^*_t(u)$ and $D_t (u)$ can be
described as follows. Let $N_+$ be the nilpotent subgroup of the
group $G$ generated by positive roots and let $H$ be Cartan subgroup.
The Borel
subgroup $B_+$ is a semi-direct product of $H$ and $N_+$.

\begin{definition}

The Poisson-Lie group $K^*_t(u)$ is the following subgroup of $B_+$:

\begin{equation}
K^*_t(u)= \{ k^* \in B_+ : k^*_{\alg{h}}=exp(i(a+iu(a))) \forsome a
\in \alg{h} \}.
\end{equation}
Here $k^*_{\alg{h}}$ is the image of the projection of $k^*$ into the
Cartan subgroup of $G$:

\begin{equation}
k^*=k^*_{\alg{h}}n \ , \ n \in \mathcal{N}_+.
\end{equation}

The Poisson structure on the group $K^*_t(u)$ is given by the
following
formula:

\begin{multline} \label{PLb}
\{ \psi , \psi' \}_{(t,u)} (k^*) = \\
= <r_{+}^t(u), \nabla_L \psi (k^*) \otimes \nabla_L \psi' (k^*) > -
<r_{+}^t(u), \nabla_R \psi (k^*) \otimes \nabla_R \psi' (k^*) > .
\end{multline}

\end{definition}
Here we embed $K^*_t(u)$ into $G$ and understand $\nabla_L$ and
$\nabla_R$ as left and right derivations on $G$.

\begin{definition}
The Heisenberg double $D_t (u)$ of the compact Poisson-Lie group $K$
is the group $G$ equipped with the following Poisson bracket:

\begin{multline} \label{DPb}
\{ \psi , \psi' \}_{(t,u)} (g) = \\
= <r_+^t(u), \nabla_L \psi (g) \otimes \nabla_L \psi' (g) > +
<r_-^t(u), \nabla_R \psi (g) \otimes \nabla_R \psi' (g) > .
\end{multline}

\end{definition}

Let us make two remarks about the Heisenberg double:

\begin{enumerate}
\item The Heisenberg double is not a Poisson-Lie group.
\item The Poisson structure \eqref{DPb} is nondegenerate.
\end{enumerate}

The group $K$ acts on $D_t(u)$ by left and right multiplications.
These  are Poisson actions.

\begin{definition}
The action  of a Poisson-Lie group $K$ on a Poisson manifold $M$ is
called a Poisson action if the map $K\times M\rightarrow M$ is a
Poisson map.
\end{definition}
Let us remark that if the Poisson bracket on the group $K$ is
nontrivial, the  Poisson action does \emph{not} preserve the Poisson
bracket on $M$.

In the case of Poisson-Lie action of the compact group $K$ the
standard notion of the moment map should be replaced by the notion of
the moment map in the sense of Lu and Weinstein \cite{Lu-Weinstein}
(Poisson-Lie moment map). The target space of this new kind of moment
map is the nonabelian Poisson-Lie group $K^*$ instead of the abelian
Poisson space $\alg{k}^*$. Let $M$ be a Poisson manifold which
carries a Poisson action of the Poisson-Lie group $K$. Define a set
of vector fields $v_X$ on $M$   for  $X\in \alg{k}$ which represent
the infinitesimal $K$-action.

\begin{definition}
A map $\mm : M\rightarrow K^*$ is called a moment map in the sense of
 Lu and Weinstein if the following equation holds true:
\begin{equation}
v_X= \mathcal{P}(\cdot , <\mm^*(\theta), X>).
\end{equation}
Here $\theta$ is the right-invariant Maurer-Cartan form on the group
$K^*$.
\end{definition}

In order to describe the moment maps for left and right
multiplication actions of $K$ on $D$ one can use the Iwasawa
decomposition:

\begin{gather}
D=G=KK^*=K^*K, \nonumber \\
g=\pi_L(g) \pi_R^*(g)= \pi^*_L(g) \pi_R(g).
\end{gather}
So, we have defined right and left projections $\pi_R$, $\pi^*_R$ and
$\pi_L$, $\pi_L^*$ from $D$ to $K$ and $K^*$.

\begin{lemma}

The moment maps for the left and right multiplication actions of the
group $K$ on $D$ are given by the projections $\pi^*_R$ and
$\pi^*_L$:

\begin{gather}
\mm^L(g) = \pi^*_L(g), \\
\mm^R(g) = (\pi^*_R(g))^{-1}.
\end{gather}

\end{lemma}

Now we turn to the Poisson structure \eqref{PLb} on $K^*$. It is
degenerate. The corresponding symplectic leaves coincide with orbits
of the dressing
action of $K$ on $K^*$ (these orbits are deformations of
coadjoint orbits in $\alg{k}^*$). One can define the dressing action
by combining the left  action of $K$ on $D$ with the right projection
from $D$ to $K^*$. More explicitly, one can introduce maps
$\rho : K \times K^* \rightarrow K$ and $\rho^* : K \times K^*
\rightarrow K^*$ defined by the following equations:

\begin{gather} \label{kkrr}
\rho(k, k^*)=\pi_R(kk^*) , \\
\rho^*(k, k^*)=\pi^*_L(kk^*). \nonumber
\end{gather}

\begin{definition}

The result of the dressing action $AD^*$ of $k \in K$ on $k^* \in
K^*$ is given
by $\rho^* ( k,k^* )$ in \eqref{kkrr}:

\begin{equation}
\AD[*]{k} k^* = \rho^* ( k,k^* ) .
\end{equation}

\end{definition}

There exists a map from $\alg{k}^*$ to $K^*$ which intertwines the
actions   $Ad^*$ and  $AD^*$. In order to define it let us
first introduce the following map:

\begin{gather}
f : K^* \rightarrow \SK
\nonumber \\
f(k^*) = k^* \overline{k^*},
\end{gather}
where $\SK$ is the following subspace in $G$:

\begin{equation}
\SK = \{ exp(iX), X \in \alg{k} \} = \{ g \in G, \overline{g}=g \} .
\end{equation}
It is easy to prove that $f$ is invertible. So, one can
define the map $e$ from $\alg{k}^*$ to $K^*$:

\begin{equation}
e(X) = f^{-1} [\exp (2it \ I^{-1}(X))].
\end{equation}
We remind that $I$ is the  map from $\alg{k}$ to $\alg{k}^*$
defined by the Killing form.

It is easy to see that $e$ intertwines the actions of $K$
on $\alg{k}^*$ and $K^*$:

\begin{equation}
\AD[*]{k} e(X) = e(\Ad[*]{k} X).
\end{equation}

The map $e$ enables us to parametrize dressing orbits
in $K^*$ by elements of the positive Weyl chamber in
Cartan subalgebra $\alg{h}$:

\begin{equation}
\orb{t}{H} =\{  k^* = \AD[*]{k} e(I(H)) , k\in K \}
\end{equation}
In fact, $\orb{t}{H}$ and $\orb{0}{H}$ are isomorphic as symplectic
spaces
\cite{Ginzburg-Weinstein}.

Now we can easily define the Poisson-Lie multiplicity
spaces. Consider a product of three dressing orbits
$\orb{t}{H_1} \times  \orb{t}{H_2} \times  \orb{t}{H_3}$.
Group $K$ acts on it in the following way

\begin{equation} \label{kOOt}
(k_1^*,k_2^*,k_3^*)
\stackrel{k}{\rightarrow}
(\AD[*]{k}k_1^*, \AD[*]{\rho(k,k_1^* )} k_2^*,
\AD[*]{\rho(\rho(k,k_1^* ),k_2^* )} k_3^* ),
\end{equation}
where $\rho$ is defined in equation \eqref{kkrr}. We refer to this
action as to diagonal dressing action.

The action \eqref{kOOt} is a Poisson one
with the Poisson-Lie moment map
$\mm (k_1^*,k_2^*,k_3^*)=k^*_1 k^*_2 k^*_3$. So one can apply the
Poisson reduction \cite{Lu} (Poisson-Lie counterpart  of the
Hamiltonian reduction). The resulting Poisson quotients are
Poisson-Lie multiplicity spaces.

\begin{definition}
Let $e_{K^*}$ be the  unit element of the group $K^*_t(u)$.
The Poisson-Lie multiplicity space $\mult{t}{H_1,H_2, H_3}$
is the quotient of  the space $\mm^{-1}(e_{K^*})$ over the action
\eqref{kOOt} of $K$:

\begin{multline}
\mult{t}{H_1,H_2,H_3} = \\
=\{ (k^*_1, k^*_2, k^*_3) \in
\orb{t}{H_1} \times \orb{t}{H_2} \times \orb{t}{H_3} :
k^*_1 k^*_2 k^*_3 = e_{K^*} \}/K.
\end{multline}

\end{definition}

The main aim of this text is to prove that there exists a global
symplectic isomorphism between $\mult{0}{H_1, H_2 ,H_3}$ and
$\mult{t}{H_1, H_2 ,H_3}$. The main tool in our proof is the moduli
space of flat connections on the sphere with three holes and a
certain real form of it.

\section{The moduli space of flat connections on the sphere with
three holes}\label{moduli}

\subsection{Atiyah-Bott description}\label{Atiyah-Bott}

Let us consider a sphere with three holes. We denote it by $\RS$ and
represent  as a complex plain without two discs $D_1$ and $D_{-1}$
around $z=1$ and $z=-1$:

\begin{equation}
\RS = \C \backslash (D_1 \cup D_{-1}).
\end{equation}

The third hole is at  infinity.

Let $\all_{\RS}$ be the space of all connections on $\RS$ with
values in the Lie algebra $\alg{g}$. $\all_{\RS}$ is a Poisson
manifold with Poisson bracket:

\begin{equation} \label{PbA}
\{ \psi ,\psi'\}=\int_\RS
<\frac{\delta \psi}{\delta A} \wedge \frac{\delta \psi'}{\delta A}>,
\end{equation}
where $ \psi(A)$ and $\psi'(A)$ are functionals on $\all_{\RS}$.
The Poisson bracket \eqref{PbA} is degenerate. In order
to describe the corresponding symplectic leaves we
fix conjugacy classes of the holonomies around
holes. Let us choose contours $\G_1$, $\G_2$ and $\G_3$ around the
points $1$, $-1$ and $\infty$ correspondingly. Define the space

\begin{equation}
\all(H_1,H_2,H_3) =
\{ A \in \all_{\RS} : \Hol(A, \G_j) \sim \exp(2it H_j)
 \},
\end{equation}
where $H_i$ belong to the positive Weyl chamber of
the Cartan subalgebra $\alg{h}$, $\Hol(A, \G_j) \sim \exp(2it H_j)$
means that $\Hol(A, \G_j)$ is conjugate to $\exp(2it H_j)$ as an
element in
$G$.

The spaces $\all(H_1,H_2 ,H_3)$ are the symplectic
leaves of  the Poisson bracket \eqref{PbA}.

The bracket \eqref{PbA} is invariant with
respect to the action of the gauge group $G_{\RS}$:

\begin{equation} \label{GS}
A \stackrel{g}{\rightarrow} g A g^{-1} -dgg^{-1}.
\end{equation}

Moreover, the gauge action is Hamiltonian, the corresponding
moment map is given by the curvature:

\begin{equation}
F(A)=dA+A^2.
\end{equation}

The symplectic quotient of the space $\all(H_1,H_2,H_3)$
is the moduli space of flat connections on the sphere with three
holes.
\begin{definition}
The moduli space of flat connections is a quotient of the space of
flat connections over the gauge action of $G_\RS$:

\begin{equation}
\moduli[]{H_1, H_2 ,H_3}=\{ A \in \all(H_1,H_2 ,H_3) :
F(A) = 0 \}
/G_\RS.
\end{equation}

\end{definition}
The space $\moduli[]{H_1,H_2 ,H_3}$ is a finite dimensional
symplectic space with Poisson structure inherited from
\eqref{PbA}.

We shall need two alternative descriptions of the Poisson
structure on the moduli space
$\moduli[]{H_1,H_2 ,H_3}$.
The first one was introduced by
Goldman \cite{Goldman}.
In order to describe the Goldman bracket let us consider the
following set of functions on the space of flat connections
on $\RS$. Denote by $\lambda$ some exact matrix  representation of
the group $G$. Then define

\begin{equation} \label{fG}
\phi_{\G}=\Tr_{\lambda} \Hol(A, \G)
\end{equation}
for a closed countor  $\G$ on the surface. Occasionally we shall
drop the sign $\lambda$ in the right hand side of this definition.
The connection $A$ being flat,  the function \eqref{fG}
depends exclusively on the homotopy class of $\G$.  Evidently,
it is invariant with respect to
the gauge transformations \eqref{GS}. So, it is a well-defined
function on the
moduli space of flat connections $\moduli[]{H_1,H_2,H_3}$
(or, more exactly, a pullback of such a function).
In the case of $G=SU(n)$ Goldman \cite{Goldman} defined the  Poisson
bracket of two such functions $\phi_{\G}$ and
$\phi_{\G^{'}}$ as follows:

\begin{equation} \label{GPB}
\{\phi_{\G}, \phi_{\G'} \}= \sum_{p \in \G \cap\G'} (-1)^{\nu(p)}
( \phi_{\G_p}-\frac{1}{n} \phi_{\G} \phi_{\G'} ).
\end{equation}
Here $\nu(p)$ is the contribution into the intersection
number of $\G$ and $\G'$
at the point $p$ and $\G_p$ is the contour produced  from
the contours $\G$
and $\G'$ by the transformation in the vicinity of
intersection point $p$ shown on Figure \ref{Goldman}.
We suppose that the contours intersect
each other transversally and in finite number of points. This is
not a restricting condition since only the homotopy class of
a contour is important.

\begin{figure}
\begin{center}
\begin{picture}(220,90)(-40,-40)

\put(-40,-40){\vector(1,1){20}}
\put(-20,-20){\line(1,1){20}}
\put(0,0){\vector(1,1){20}}
\put(20,20){\line(1,1){20}}
\put(40,-40){\vector(-1,1){20}}
\put(20,-20){\line(-1,1){20}}
\put(0,0){\vector(-1,1){20}}
\put(-20,20){\line(-1,1){20}}

\put(100,-40){\line(1,1){20}}
\put(120,-20){\vector(0,1){20}}
\put(120,0){\line(0,1){20}}
\put(120,20){\line(-1,1){20}}
\put(180,-40){\line(-1,1){20}}
\put(160,-20){\vector(0,1){20}}
\put(160,0){\line(0,1){20}}
\put(160,20){\line(1,1){20}}

\put(0,0){\circle*{3}}
\put(140,0){\circle*{3}}

\put(60,-6){$\rightarrow$}

\put(-3,-14){$p$}
\put(137,-14){$p$}
\put(-27,-40){$\G$}
\put(-27,32){$\G'$}
\put(17,-40){$\G'$}
\put(17,32){$\G$}
\put(106,-4){$\G_p$}
\put(163,-4){$\G_p$}

\end{picture}
\end{center}
\caption{} \label{Goldman}
\end{figure}

One can easily derive bracket \eqref{GPB} from \eqref{PbA}.

More complicated functions on the moduli space can be defined by
replacing countors by arbitrary graphs. The Poisson brackets of these
functions similar to the Goldman's brackets can be found in
\cite{AMR}.

\subsection{Fock-Rosly theorem}

Another description of the Poisson bracket on the moduli space
of flat connections was introduced by Fock and Rosly in
\cite{Fock-Rosly}. They used a notion of a ciliated fat graph
which we are going to review now. Let us draw an oriented graph $L$
on the  surface $\RS$  such that each face of
the graph is homeomorphic either to a disk or to a disk with one
hole.
For instance, one can consider  a triangulation of the surface.
The graph $L$ has a cyclic order of ends of edges at each
vertex which is
inherited from the orientation of the surface. Such a graph
is called a fat graph. It is convenient to introduce notations $s(l)$
and $t(l)$ for the starting and end points of the oriented edge $l$.

Let $\all_L$ be the space of graph connections. A graph
connection is  an assignment of an
element $a_l$ of the group $G$ to each edge $l$ of the graph.
Correspondingly, an element of the graph gauge group $G_L$
is an assignment of an
element $g_{\alpha} \in G$ to each vertex $\alpha$ of the graph.
The graph gauge  group naturally acts  on the space of graph
connection:

\begin{equation}
a_l \longrightarrow g_{t(l)} a_l g_{s(l)}^{-1}.
\end{equation}

We  call a graph connection flat if its holonomy  around each
empty face is equal to identity in $G$ and the  holonomy around the
face
containing the hole $\G_i$ belongs to the same conjugacy class as
$exp(2i tH_i)$. We denote the space of flat graph connections by
$\allflat_L(H_1, H_2, H_3)$. It is a standard fact
that the quotient of the space of flat graph connections over the
graph gauge group is equal to the moduli space
of flat connections on $\RS$:

\begin{equation}
\moduli[]{H_1, H_2, H_3}=\allflat_L(H_1, H_2, H_3)/G_L.
\end{equation}

Fock and Rosly \cite{Fock-Rosly} introduced a   Poisson structure on
the space $\all_L$.  This structure depends on a linear order of the
ends of edges at each vertex of the graph.
One can fix the linear  order by putting
a small cilium at the vertex and choosing the edge just after the
cilium
to be the first in the enumeration (see Figure \ref{cilium}).

\begin{figure}
\begin{center}
\begin{picture}(80,95)(-30,-40)

\put(0,0){\vector(0,1){22}}
\put(0,22){\line(0,1){22}}
\put(0,0){\line(2,1){20}}
\put(40,20){\vector(-2,-1){20}}
\put(0,0){\line(1,-2){10}}
\put(20,-40){\vector(-1,2){10}}
\put(0,0){\vector(-1,-1){15}}
\put(-15,-15){\line(-1,-1){15}}

\put(0,0){\circle*{3}}
\put(0,0){\line(-1,1){10}}

\put(3,35){$l_1$}
\put(35,7){$l_2$}
\put(20,-36){$l_3$}
\put(-23,-33){$l_4$}

\end{picture}
\end{center}
\caption{} \label{cilium}
\end{figure}

A fat graph with this additional structure is called a
ciliated fat graph.

We need some notations.
Let $N(L)$ be the set of all vertices of the graph $L$ and $E(L)$ be
the set of ends of  edges of $L$. We  denote by $E(n)$ the subset of
$E(L)$ which consists of the ends of edges incident to a given vertex
$n\in N(L)$.

Let us introduce vector fields
\begin{equation}
X_{\alpha}=
\begin{cases}
\nabla^L_l & \text{for $\alpha$ being the target end of the edge $l$}
\\
-\nabla^R_l & \text{for $\alpha$ being the starting end of the edge
$l$}.
\end{cases}
\end{equation}
Here
$\nabla^L_l$ and $\nabla^R_l$ are left- and
right-invariant derivations on the copy of the group $G$ assigned
to an edge $l \in L$.

Now we are ready to write down the Poisson bracket of
two functions $\psi$ and $\psi'$ on the space of graph
connection. It is given by the following formula:

\begin{multline}\label{FRPb}
\{\psi , \psi' \}= \\ =
\sum_{n \in N(L)} ( \sum_{\substack{\alpha, \beta \in E(n) \\
\alpha <\beta}} <r^t_+,X^1_{\alpha} \psi \wedge X^2_{\beta}\psi'>+
\frac{1}{2}\sum_{\alpha \in E(n)}
<r^t_+,X^1_{\alpha}\psi \wedge X^2_{\alpha}\psi'>),
\end{multline}
where $r^t_+$ is the classical $r$-matrix (\ref{r+t}).
Here we use tensor notations $X^1=X \otimes 1$ and
$X^2=1 \otimes X$. We write $\alpha < \beta$ in the
sense of the linear order in $E(n)$.

The main result of \cite{Fock-Rosly} is described by the following
theorem.
 \begin{theorem} \

\begin{enumerate}
\item
The Poisson bracket \eqref{FRPb} induces a Poisson bracket on the
quotient space $\all_L/G_L$.
\item
The moduli space of flat connections is embedded into the quotient
space
$\all_L/G_L$ as a symplectic leaf so that the induced Poisson bracket
coincides with the Goldman bracket.
\end{enumerate}
\end{theorem}

\subsection{The hyperbolic moduli space}

Let $\tau$ be a second order anti-automorphism
(in the sense of complex structure) of the Riemann
surface $\RS$:

\begin{equation}
\tau : \RS \longrightarrow \RS \ , \ \tau^2=id.
\end{equation}
Combining $\tau$ with the   anti-involution of $\alg{g}$ one can
define an anti-automorphism $\aut$ on the space of connections
$\all_{\RS}$:

\begin{equation} \label{defaut}
\aut(A)(z) = - \overline{A (\tau(z))}.
\end{equation}

One of the main objects of our interest in this
text is the $\aut$-invariant  subspace
in $\all_{\RS}$:

\begin{equation}\label{As}
\all_{\RS}^{\aut}=\{ A \in \all_{\RS}: \aut(A)=A \}.
\end{equation}
The space $\all_{\RS}^{\aut}$ is acted on by
the subgroup $G_{\RS}^{\aut}$ of gauge group
$G_{\RS}$:

\begin{equation}\label{Gs}
G_{\RS}^{\aut}=\{ g \in G_{\RS}: \overline{g(z)}=g^{-1}(\tau(z)) \}.
\end{equation}

If the surface $\RS$ has holes we demand that holonomies
around the holes belong to some fixed conjugacy classes (see
subsection \ref{Atiyah-Bott}).  We assume that the map $\tau$ is an
anti-automorphism of the Riemann surface with holes. We shall suppose
that each particular hole is preserved by $\tau$. Choose  the
contours $\G_i$ which are used to compute holonomies to be
$\tau$-invariant.  We can define the following subspace in the space
$\all_{\RS}^{\aut}$:

\begin{equation} \label{hyp}
\all^{\aut}_{\RS}(H_1, H_2 ,H_3) =
\{ A \in \all_{\RS}^{\aut} : \Hol(\G_j, A) \sim \exp (2itH_j)
, H_j\in W_+ \}.
\end{equation}
We also define the subspace of flat connections in
$\all^{\aut}_{\RS}(H_1, H_2 ,H_3)$

\begin{equation}
\allflat^{\aut}_\RS(H_1, H_2 ,H_3)=
\{ A \in \all^{\aut}(H_1, H_2 , H_3):
F(A)=0 \}.
\end{equation}
and the \emph{hyperbolic} moduli space of flat connections.

\begin{definition}
The hyperbolic moduli space of flat connections is the quotient
of the space $\allflat^{\aut}_\RS(H_1, H_2 , H_3)$ over the action of
the gauge group \eqref{Gs}
\begin{equation}
\moduli{H_1, H_2 , H_3}=
\allflat^{\aut}_\RS(H_1, H_2 , H_3)/G^{\aut}_\RS.
\end{equation}
\end{definition}
It is obvious that Poisson bracket \eqref{PbA}
can be restricted to the space $\all^{\aut}$ and
so the hyperbolic moduli space is
a symplectic quotient of $\all^{\aut}(H_1, H_2 , H_3)$.

\begin{remark}
It is instructive to compare the hyperbolic moduli space of flat
connections with the moduli space of flat connections with values in
the compact form $\alg{k}$.
In the latter case one considers the space of connections restricted
by the condition
\begin{equation}
A(z) = - \overline{A (z)}.
\end{equation}
The compact gauge group is defined as
\begin{equation}
G_{\RS}^{\alg{k}}=\{ g \in G_{\RS}: \overline{g(z)}=g^{-1}(z)\}.
\end{equation}
Finally,  the holonomies around holes are restricted to some
conjugacy classes in the compact group $K$:
\begin{equation}
\all^{\alg{k}}_{\RS}(H_1, H_2 ,H_3) =
\{ A \in \all_{\RS}^{\alg{k}} : \Hol(\G_j, A)\sim \exp (2tH_j)
, H_j\in W_+ \}.
\end{equation}
The eigenvalues of the holonomy matrices in this case are unitary
numbers. For this reason one can call the moduli space
$\moduli[\alg{k}]{H_1, H_2 , H_3}$  \emph{elliptic} moduli space of
flat connections. On the contrary, the eigenvalues of the holonomy
matrices defined by \eqref{hyp} are positive real numbers. They
define a set of hyperbolic elements of $G$ and we call the
corresponding moduli space the \emph{hyperbolic} moduli space of flat
connections.
\end{remark}

Now we have to adjust the Fock and Rosly
description to the hyperbolic moduli space. We choose a ciliated fat
graph $L$
to be invariant with respect to the anti-automorphism $\tau$:

\begin{equation}
\tau(L)=L.
\end{equation}

The analogue of the condition (\ref{As})
for graph connections is

\begin{equation}\label{as}
a_l=(\overline{a}_{\tau(l)})^{ -1},
\end{equation}
where $\tau(l)$ is the image
of the edge $l$ under the action of $\tau$.

Similarly, the  condition for the graph gauge group
is the following:

\begin{equation}\label{gs}
g_{\alpha}=(\overline{g}_{\tau (\alpha)})^{ -1}
\end{equation}
for any vertex $\alpha$ and its image $\tau(\alpha)$

We leave it as an exercise for the reader to check that
the Fock-Rosly bracket (\ref{FRPb}) can be
restricted to graph connections obeying (\ref{as}).
In proving this fact it is
important to note that the linear order  at a
vertex reverses under the action of $\tau$ and to use the identity

\begin{equation}
r^t_-=-P(r^t_+)=\overline{r^t_+}.
\end{equation}

So, one can obtain the hyperbolic moduli space
$\moduli{H_1, H_2 ,H_3}$ as a
quotient of the space of flat graph connections obeying condition
(\ref{as})
over the subgroup (\ref{gs}) of the graph gauge group.

In the case of a sphere with three holes $\RS = \C \backslash (D_1
\cup D_{-1})$ we choose the anti-automorphism $\tau$ to be
represented by the complex conjugation:
\begin{equation}
\tau(z)=\overline{z}.
\end{equation}
The holes are preserved by the action of $\tau$. The corresponding
hyperbolic moduli space is denoted by $\moduli{H_1, H_2 ,H_3}$.

\section{Equivalence of Poisson-Lie and compact multiplicity
spaces}\label{isomorphisms}

A triple  $(H_1, H_2, H_3)$ of elements of the positive Weyl chamber
of $\alg{h}$ defines three different symplectic spaces:
\begin{enumerate}
\item the symplectic multiplicity space $\mult{0}{H_1, H_2 ,H_3}$,
\item  the Poisson-Lie multiplicity space $\mult{t}{H_1, H_2 ,H_3}$,
\item the hyperbolic moduli space of flat connections  $\moduli{H_1,
H_2 ,H_3}$.
\end{enumerate}

The goal of this paper is to compare these three manifolds as
symplectic spaces.

\subsection{Small values of $H_1$, $H_2$ and $H_3$}
In this subsection we restrict ourselves to the situation when the
elements $H_1$, $H_2$ and $H_3$ are sufficiently close to $0$ in
$\alg{h}$.

\begin{theorem}

Let $t$ be a real number. There exists a neighborhood $U$ of $0$ in
$\alg{h}$ such that
\begin{equation}
\mult{0}{H_1, H_2 ,H_3}\approx  \mult{t}{H_1, H_2 ,H_3}\approx
\moduli{H_1, H_2 ,H_3}
\end{equation}
if  $H_1, H_2$ and $H_3$ belong to $U$.
\end{theorem}

\begin{proof}

In \cite{Jeffrey} it has been shown that $\mult{0}{H_1, H_2 ,H_3}$ is
isomorphic to the elliptic moduli space $\moduli[\alg{k}]{H_1, H_2
,H_3}$ if $H_1, H_2$ and $H_3$ are sufficiently small. The proof for
the hyperbolic moduli space follows the same scheme.

Here we prove only the isomorphism of  $\mult{0}{H_1, H_2 ,H_3}$ and
$\mult{t}{H_1, H_2 ,H_3}$. The idea of the proof is very similar to
the one in \cite{Jeffrey} (see Theorem 6.6).
Consider the direct product of three copies of the Heisenberg double
$D^3$. Define the action of the group $K^4$ on $D^3$ in the following
fashion:
\begin{equation} \label{act}
(d_1, d_2, d_3)\rightarrow (k_0^{-1}d_1k_1, k_0^{-1}d_2k_2,
k_0^{-1}d_3k_3).
\end{equation}
This is a Poisson action. The Poisson-Lie moment map for this action
is given by
\begin{equation}
\mm: (d_1, d_2, d_3)\rightarrow (\mm^L(d_1) \mm^L(d_2) \mm^L(d_3),
\mm^R(d_1), \mm^R(d_2), \mm^R(d_3)).
\end{equation}
The reduced space corresponding to the level $(e_{K^*}, \exp(2itH_1),
\exp(2it H_2), \exp(2it H_3) )$ of the moment map coincides with
$\mult{t}{H_1, H_2 ,H_3}$.

It is known that for any Poisson action of a compact Poisson-Lie
group $K$ on a symplectic manifold which possesses an equivariant
moment map
one can introduce a $K$-invariant symplectic structure on the same
manifold  such that the reduced spaces do not change \cite{AA}. Let
us apply this construction to $D^3$ equipped with the $K^4$-action
\eqref{act}. Then by applying the theorem of  Lerman and Sjamaar
\cite{Lerman-Sjamaar} (see Proposition 2.5) one can prove that there
exists a small neighborhood of the zero level set of the moment map
in $D^3$ which is isomorphic to a small neighborhood of the zero
section in $T^*K^3$. The corresponding symplectic quotients of
$T^*K^3$ coincide with symplectic multiplicity spaces $\mult{0}{H_1,
H_2, H_3}$.

Thus we conclude that
\begin{equation}
\mult{t}{H_1, H_2 ,H_3}\approx \mult{0}{H_1, H_2, H_3}
\end{equation}
for $H_1$, $H_2$ and $H_3$ being sufficiently close to $0$.

\end{proof}

\subsection{The Poisson map from the multiplicity space to the moduli
space}\label{MultToModuli}

Let us consider the map $\xi$ from the
multiplicity space $\mult{0}{H_1,H_2,H_3}$
to the hyperbolic moduli space
$\moduli{H_1,H_2 ,H_3}$ given by
the following formula:

\begin{gather}
\xi :
\mult{0}{H_1,H_2,H_3} \rightarrow
\moduli{H_1,H_2 ,H_3}
\nonumber \\
\xi (X_1, X_2, X_3) = A(z) dz,
\\ \intertext{where}
\label{AXXX}
A(z) = \frac{X_1}{z-1} +
\frac{X_2}{z+1}.
\end{gather}
The map $\xi$ intertwines the diagonal coadjoint action of $K$ on the
triple $(X_1, X_2, X_3)$ and the action of the constant gauge
transformations on $A(z)$.  The connection $A(z)$ is flat and it
satisfies the condition $\overline{A(z)}=-A(\overline{z})$. Thus,
$\xi$ defines a map from the multiplicity space to the hyperbolic
moduli space of flat connections.

\begin{theorem}

The map $\xi$ is a Poisson map.

\end{theorem}

Before proving the theorem let us prove the following lemma.

\begin{lemma}
The Poisson bracket on
$\moduli{H_1,H_2 ,H_3}$
induced from
$\mult{0}{H_1,H_2,H_3}$ by the map $\xi$
is given by the following formula:

\begin{equation} \label{lem}
\{ \psi ,\psi'\}=\int_{\RS \times \RS} dz \ dw
<[\frac{\delta \psi}{\delta A(z)} , \frac{\delta \psi'}{\delta
A(w)}],
\frac{A(z) - A(w)}{z - w}>.
\end{equation}

\end{lemma}

\begin{remark}
The Poisson bracket \eqref{lem} is of the Sklyanin type
\cite{Sklyanin}. In the connection with the moduli spaces it has been
proposed in \cite{Fock-Rosly}.
\end{remark}

\begin{proof}[Proof of Lemma]

The proof is a straightforward calculation using
formula \eqref{KKPb} for Poisson bracket on $\alg{g}^*$
and formula \eqref{AXXX} for $\xi$:

\begin{multline}
\{ \psi ,\psi'\}=
\int_{\RS \times \RS} dz \ dw
(<[\frac{\delta \psi}{\delta A(z)} \frac{\partial A(z)}{\partial
X_1},
\frac{\delta \psi'}{\delta A(w)} \frac{\partial A(w)}{\partial X_1}],
X_1> + \\
+ <[\frac{\delta \psi}{\delta A(z)} \frac{\partial A(z)}{\partial
X_2},
\frac{\delta \psi'}{\delta A(w)} \frac{\partial A(w)}{\partial X_2}],
X_2>)=
\end{multline}

\begin{multline}
=\int_{\RS \times \RS} dz \ dw
(<[\frac{\delta \psi}{\delta A(z)}, \frac{\delta \psi'}{\delta
A(w)}],
\frac{X_1}{(z-1)(w-1)}> + \\
+ <[\frac{\delta \psi}{\delta A(z)}, \frac{\delta \psi'}{\delta
A(w)}],
\frac{X_2}{(z+1)(w+1)}>)=
\end{multline}

\begin{equation}
=\int_{\RS \times \RS} dz \ dw
<[\frac{\delta \psi}{\delta A(z)} , \frac{\delta \psi'}{\delta
A(w)}],
\frac{A(z) - A(w)}{z - w}>
\end{equation}

\end{proof}

\begin{proof}[Proof of the theorem]

In order to prove the theorem we calculate the Poisson
bracket of two Goldman functions using \eqref{lem} and show that
it coincides with the formula \eqref{GPB}.
Let $\phi_{\G}$ and $\phi_{\G'}$ be two Goldman functions:

\begin{gather}
\phi_{\G}=\Tr \Hol(A, \G, t_{0}),
\\
\phi_{\G'}=\Tr \Hol(A, \G', s_0).
\end{gather}
We parametrize the closed contours $\G$ and $\G'$ by
$t$ and $s$ correspondingly and denote the starting
point for the holonomies by $t_0$ and $s_0$.

Now we are ready to calculate the Poisson bracket of
$\phi_{\G}$ and $\phi_{\G'}$.

\begin{equation}
\{ \phi_{\G}, \phi_{\G'} \} =
\int_{\RS \times \RS} dz \ dw
<[\frac{\delta \phi_{\G}}{\delta A(z)} ,
\frac{\delta \phi_{\G'}}{\delta A(w)}],
\frac{A(z) - A(w)}{z - w}>=
\end{equation}

\begin{multline}
=\oint \oint dt ds
\frac{\partial z(t)}{\partial t} \frac{\partial z(s)}{\partial s}
\Tr_{12} \ (\Hol(A, \G, t_0)^1 \Hol(A, \G', s_0)^2  \frac{[C^{12},
A^1(z(t))]}{z(t) - z(s)}+ \\
+
\Hol(A, \G, t_0)^1 \Hol(A, \G', s_0)^2 \frac{[C^{12},
A^2(z(s))]}{z(t) - z(s)})=
\end{multline}

\begin{equation} \label{dtds}
= \oint \oint dt ds
\frac{\partial z(t)}{\partial t} \frac{\partial z(s)}{\partial s}
\frac{(\frac{\partial}{\partial z(t)} + \frac{\partial}{\partial
z(s)})
\Tr_{12}\ (C^{12} \Hol(A, \G, t)^1 \Hol(A, \G', s)^2)}{z(t) - z(s)}
\end{equation}

Here we use the following notations: $X^1=X \otimes 1$,
$X^2=1 \otimes X$ for $X \in \alg{g} \text{ or } G$,
$\Tr_{12}$ is the trace in the tensor product of two spaces
and $C^{12}=\sum_a t^a \otimes t^a$ is the tensor Casimir operator.

In order to treat integral \eqref{dtds} let us enumerate the
intersection points of $\G$ and $\G'$:

\begin{equation}
\{ p_n \} = \G \cap \G' .
\end{equation}
We suppose that the contours intersect transversally in finite number
of points.

Let $U_n$ and $V_n$ be small open neighborhoods of the points
$p_n$ such that $\overline{U_n} \subset V_n$ and let $\eta_n$ be
smooth functions on $\RS$ such that

\begin{equation}
\eta_n \mid_{U_n} = 1, \eta_n \mid_{\RS \backslash V_n} = 0.
\end{equation}

Then

\begin{equation}
\oint \oint dt ds
\frac{\partial z(t)}{\partial t} \frac{\partial z(s)}{\partial s}
\frac{( \frac{\partial}{\partial z(t)} + \frac{\partial}{\partial
z(s)})
\Tr_{12}\ (C^{12} \Hol(A, \G, t)^1 \Hol(A, \G', s)^2)}{z(t) - z(s)}=
\end{equation}

\begin{multline}
=\oint \oint dt ds
\frac{\partial z(t)}{\partial t} \frac{\partial z(s)}{\partial s}
\cdot \\ \cdot
\frac{( \frac{\partial}{\partial z(t)} + \frac{\partial}{\partial
z(s)})
(\sum_n \eta_n (z(t)) \eta_n (z(s))
\Tr_{12}\ (C^{12} \Hol(A, \G, t)^1 \Hol(A, \G', s)^2)) }{z(t) -
z(s)}+
\end{multline}
\begin{multline*}
+ \oint \oint dt ds
\frac{\partial z(t)}{\partial t} \frac{\partial z(s)}{\partial s}
\cdot \\ \cdot
\frac{( \frac{\partial}{\partial z(t)} + \frac{\partial}{\partial
z(s)})
((1-\sum_n \eta_n (z(t)) \eta_n (z(s)))
\Tr_{12}\ (C^{12} \Hol(A, \G, t)^1 \Hol(A, \G', s)^2)) }{z(t) -
z(s)}=
\end{multline*}

\begin{equation}
=J_1 + J_2
\end{equation}

Let us first consider the integral $J_2$. The integrand
has no singularities. So one can integrate by parts:

\begin{multline}
J_2 = \oint \oint dt ds
\frac{\partial z(t)}{\partial t} \frac{\partial z(s)}{\partial s}
\cdot \\ \cdot
\frac{( \frac{\partial}{\partial z(t)} + \frac{\partial}{\partial
z(s)})
((1-\sum_n \eta_n (z(t)) \eta_n (z(s))) \Tr_{12} \ (C^{12} \Hol(A,
\G, t)^1 \Hol(A, \G', s)^2 ) ) }{z(t) - z(s)}=
\end{multline}

\begin{multline}
= \oint \oint dt ds
\frac{\partial z(t)}{\partial t} \frac{\partial z(s)}{\partial s}
((1-\sum_n \eta_n (z(t)) \eta_n (z(s)))
\cdot \\ \cdot
\Tr_{12} \ (C^{12} \Hol(A, \G, t)^1 \Hol(A, \G', s)^2))
(\frac{\partial}{\partial z(t)} +\frac{\partial}{\partial z(s)})
\frac{1}{z(t) - z(s)} =0
\end{multline}

To calculate $J_1$ let us choose local coordinates $z_n$ in the
vicinities $V_n$ in such a way that  $z_n=t - is$.
Depending on the intersection number of the contours at $p_n$
this can require the change of the orientation
and so the change of the overall sign. Now the integral $J_1$ takes
the
following form:

\begin{multline}
J_1 = \oint \oint dt ds
\frac{\partial z(t)}{\partial t} \frac{\partial z(s)}{\partial s}
\cdot \\ \cdot
\frac{( \frac{\partial}{\partial z(t)} + \frac{\partial}{\partial
z(s)})
( \sum_n  \eta_n (z(t)) \eta_n (z(s))
\Tr_{12} \  (C^{12} \Hol(A, \G, t)^1 \Hol(A, \G', s)^2 )) }{z(t) -
z(s)}=
\end{multline}

\begin{multline}
= \sum_n (-1)^{\nu(p_n)}
\int_{V_n} dz_n d\overline{z}_n \cdot \\ \cdot
\frac{ \frac{\partial}{\partial \overline{z_n} }
( \eta_n (\re z_n) \eta_n (\im z_n)
\Tr_{12} \ (C^{12} \Hol(A, \G, \re z_n)^1 \Hol(A, \G', \im z_n)^2 ) )
}{z_n}=
\end{multline}

\begin{equation} \label{ans}
= \sum_{p \in \G \cap\G'} (-1)^{\nu(p)}
\Tr_{12} \ (C^{12} \Hol(A, \G, p)^1 \Hol(A, \G', p)^2 ).
\end{equation}

For $K=SU(n)$ the right hand side can be rewritten as
\begin{equation} \label{EOP}
 \sum_{p \in \G \cap\G'} (-1)^{\nu(p)}
(\phi_{\G_p}-\frac{1}{n}\phi_{\G_1}\phi_{\G_2}).
\end{equation}
Here $\nu(p)$ is the intersection number at the point $p$ and
the contour $\G_p$ is described in section \ref{moduli} (see
Figure \ref{Goldman}). Formula \eqref{EOP} coincides with
the Goldman definition \eqref{GPB}. So the map $\xi$ is
a Poisson map.

\begin{remark}
Formula \eqref{ans} is valid for any simple Lie algebra $\alg{g}$. It
always coincides with the Goldman bracket. In fact, one can repeat
the same proof for an arbitrary function on the moduli space of flat
connections \cite{Fock-Rosly}. For the sake of simplicity we present
only the calculation for the Goldman functions in the case of
$K=SU(n)$.
\end{remark}

\end{proof}

\subsection{The Poisson map from the moduli space to the Poisson-Lie
multiplicity space}\label{ModuliToMult}

In this subsection we describe a Poisson map from the hyperbolic
moduli space of flat connections on the sphere with three holes into
the Poisson-Lie multiplicity space.

Let us draw a graph on the complex plain as shown at Figure
\ref{graph}.  It is convenient to enumerate the edges of this graph
as $\G_i$, $\overline{\G}_i, i=1, 2, 3$ and the vertices as  $P_i,
i=0, 1, 2, 3$. Given a flat connection satisfying the condition
$\overline{A(z)}=-A(\overline{z})$, one can define a triple $(g_1,
g_2, g_3)\in G^3$ as
\begin{equation}
\CH : A(z)\rightarrow \{ g_i \} \ , \ g_i= \Hol(A, \G_i)
\end{equation}

The space $G^3$ is acted on by the group $K^4$ in the following
fashion:
\begin{equation} \label{K4}
g_i\rightarrow k_{i-1}^{-1} g_i k_{i}.
\end{equation}
The map $\CH$ intertwines the action of the gauge group \eqref{GS}
and the action \eqref{K4} understood as an action of the gauge group
by means of the projection:
\begin{gather}
\pi_4: G_{\RS}\rightarrow K^4 \nonumber \\
\pi_4: g(z)\rightarrow \{ g(P_i)\}.
\end{gather}
So, the moduli space is mapped into the quotient $G^3/K^4$.

Let us split the projection $G^3\rightarrow G^3/K^4$ into two pieces
corresponding to the subgroups $K_0\simeq K$ and $K_{123}\simeq K^3$
in $K^4$:
\begin{gather}
K_0=\{ (k, e_K, e_K, e_K)\in K^4 \}, \\
K_{123}=\{ (e_K, k_1, k_2,  k_3)\in K^4 \}. \nonumber
\end{gather}
As the group $G$ is the double for $K$, the quotient $G^3/K^3\simeq
(K^*)^3$. More explicitly:
\begin{gather}
(g_1, g_2, g_3)\rightarrow (k_1^*, k_2^*, k_3^*) \nonumber \\
k^*_1=\pi^*_L (g_1),   \\
k^*_2=\pi_L^*(\pi_R(g_1) g_2), \nonumber \\
k^*_3 =\pi_L^*(\pi_R(\pi_R(g_1) g_2) g_3). \nonumber
\end{gather}
Here $\pi_L^*$ and $\pi_R$ were defined in subsection \ref{PLmult}.
It is easy to see that $K_0$ acts on the quotient space $(K^*)^3$ by
the diagonal dressing action \eqref{kOOt}. In this way we define the
map
\begin{equation}
\chi: \moduli{H_1, H_2, H_3} \rightarrow (K^*)^3/K.
\end{equation}
Observe that  the image of the map $\chi$ satisfies two additional
conditions:
\begin{enumerate}
\item The connection $A(z)$ being flat, the product of $k^*_i$ is
equal to
identity:
\begin{equation}
k^*_1 k^*_2 k^*_3=e_{K^*},
\end{equation}
\item As holonomies around the holes are restricted to certain
conjugacy classes in $G$ then
\begin{equation}
k^*_i \in \orb{t}{H_i} .
\end{equation}
\end{enumerate}
We conclude that  $\chi$ actually maps the hyperbolic moduli space
into the Poisson-Lie multiplicity space.

\begin{figure}
\begin{center}
\begin{picture}(270,160)(-100,-80)

\multiput(-40,0)(80,0){2}{\oval(80,40)}
\put(20,0){\oval(200,120)}

\multiput(-40,20)(81,0){2}{\vector(1,0){0}}
\multiput(-40,-20)(81,0){2}{\vector(1,0){0}}
\multiput(15,-60)(0,120){2}{\vector(-1,0){0}}

\put(-100,0){\vector(1,0){260}}

\multiput(-80,0)(40,0){6}{\circle*{3}}

\put(-47,24){$\G_2$}
\put(34,24){$\G_1$}
\put(-47,-33){$\overline{\G}_2$}
\put(34,-33){$\overline{\G}_1$}
\put(13,-73){$\overline{\G}_3$}
\put(13,64){$\G_3$}

\put(-77,3){$P_2$}
\put(2,3){$P_1$}
\put(82,3){$P_0$}
\put(122,3){$P_3$}

\multiput(-80,0)(80,0){3}{\line(2,-1){10}}
\put(120,0){\line(2,-1){10}}
\put(-50,-11){$-1$}
\put(38,-11){$1$}
\put(-8,-10){$0$}

\put(145,-13){$\re z$}

\end{picture}
\end{center}
\caption{} \label{graph}
\end{figure}

\begin{lemma}
The map
$ \chi: \moduli{H_1, H_2, H_3}\rightarrow \mult{t}{H_1, H_2, H_3}$
is a Poisson map.
\end{lemma}

\begin{proof}
Let us consider the space $G^3$ as a real manifold. The ciliated
graph (Figure \ref{graph}) defines a Fock-Rosly bracket on $G^3$.
This bracket descends to $G^3/K^4$. One of symplectic leaves in
$G^3/K^4$ coincides with the image of the subspace
\begin{equation} \label{con1}
\allflat=
\{ (g_1, g_2, g_3)\in G^3 \ , \ g_1g_2g_3=e_G \ , \ g_j\overline{g_j}
\sim \exp{2itH_j} \}
\end{equation}
and is isomorphic to the hyperbolic moduli space $\moduli{H_1, H_2,
H_3}$.

On the other hand, one can check  that the projection $G^3\rightarrow
(K^*)^3\simeq G^3/K_{123}$ is a Poisson map if one takes the
Fock-Rosly bracket on $G^3$ and the direct sum of canonical brackets
on three copies of $(K^*)^3$.
By definition, the Poisson-Lie multiplicity space $\mult{t}{H_1, H_2,
H_3}$ is a symplectic leaf  in $(K^*)^3/K$ which is singled out by
the following conditions:
\begin{gather}
  \label{con2}
 k^*_1k^*_2k^*_3=e_{K^*} , \\
k_i \in \orb{t}{H_i}. \label{con3}
\end{gather}

As  $G^3/K^3=(K^*)^3$ as Poisson spaces and conditions \eqref{con1}
coincide with conditions \eqref{con2} and \eqref{con3}, we conclude
that
\begin{equation}
\moduli{H_1, H_2, H_3}\approx \mult{t}{H_1, H_2, H_3}.
\end{equation}

\end{proof}

Therefore we arrive at the following theorem.

\begin{theorem}
The hyperbolic moduli space $\moduli{H_1, H_2, H_3}$ is isomorphic to
the Poisson-Lie multiplicity space $\mult{t}{H_1, H_2, H_3}$ as a
symplectic space.
\end{theorem}

\subsection{The main theorem}

Let us summarize all the information about the hyperbolic moduli
space and the multiplicity spaces which we collected by now. We
consider three families of smooth symplectic manifolds $\mult{0}{H_1,
H_2, H_3}$, $\mult{t}{H_1, H_2, H_3}$ and $\moduli{H_1, H_2, H_3}$
and the following maps:
\begin{enumerate}
\item The Poisson map $\xi: \mult{0}{H_1, H_2, H_3}\rightarrow
\moduli{H_1, H_2, H_3}$,
\item The symplectic isomorphism $\chi: \moduli{H_1, H_2,
H_3}\rightarrow \mult{t}{H_1, H_2, H_3}$,
\item The symplectic isomorphisms between all three spaces for small
values of $H_1$, $H_2$ and $H_3$. These isomorphisms are known to
exist but are not known explicitly.
\end{enumerate}

It is our goal to prove that the map $\xi$ is also a symplectic
isomorphism.

\begin{lemma}
The map $\xi$ is a symplectic isomorphism.
\end{lemma}

\begin{proof}
Both spaces $\moduli{H_1, H_2, H_3}$ and $\mult{0}{H_1, H_2, H_3}$
have the same dimension equal to $(\Dim K - 3\rank K)$.
Thus, the map $\xi$ is a Poisson map between symplectic spaces of the
same dimension. Hence, it defines a covering of the space
$\moduli{H_1, H_2, H_3}$ by the space $\mult{0}{H_1, H_2, H_3}$. We
know that for small values of $H_1$,  $H_2$ and $H_3$ these spaces
are actually isomorphic as symplectic spaces. Both $\moduli{H_1, H_2,
H_3}$ and $\mult{0}{H_1, H_2, H_3}$ depend continuously on $H_1$,
$H_2$ and $H_3$. Therefore, the degree of the covering can not
change.  We conclude that the map $\xi$ is invertible.
\end{proof}

This completes the proof of the main result of this paper.

\begin{theorem}
The following three symplectic spaces
\begin{enumerate}
\item  the multiplicity space $\mult{0}{H_1, H_2, H_3}$,
\item  the Poisson-Lie multiplicity space $\mult{t}{H_1, H_2, H_3}$,
\item  the hyperbolic moduli space of flat connections on the sphere
with three holes $\moduli{H_1, H_2, H_3}$
\end{enumerate}
are isomorphic as symplectic manifolds.
\end{theorem}

Let us remark that both isomorphisms $\xi$ and $\chi$  which play the
role  in the last theorem are defined by explicit formulae (see
previous subsections).

\subsection*{Acknowledgements}
A.A. is grateful to V.Rubtsov and J.Teschner for useful discussions.
A.M. thanks  Prof. A.J. Niemi for hospitality at Uppsala.


\begin{thebibliography}{99}
\bibitem{AA} A.Yu. Alekseev:  On Poisson actions of compact Lie
groups on symplectic manifolds. \mbox{dg-ga/9602001}


\bibitem{AMR} J. Andersen, J.Mattes, N.Reshetikhin:  Poisson
structure on the moduli space of flat connections and chord diagrams.
Berkeley preprint


\bibitem{Arnold}
V.I. Arnold:
Mathematical methods of classical mechanics.
Springer-Verlag 1980

\bibitem{Fock-Rosly}
V.V. Fock, A.A. Rosly:
Poisson Structure on the moduli of flat
connections on Riemann surfaces and
$r$-matrix.
preprint ITEP 92-72 (1992);

Flat connections and polyubles.
Teor. Mat. Fiz. \textbf{95} (1993) 228


\bibitem{Ginzburg-Weinstein}
V.L.Ginzburg, A.Weinstein:
 Lie Poisson structure on some Poisson Lie groups. J. Amer. Math.
Soc. \textbf{5} (1992) 445

\bibitem{Goldman}
W. Goldman:
Invariant  functions on Lie groups
and Hamiltonian flows of surfaces group
representations. Invent. Math. \textbf{85} (1986) 263

\bibitem{Hitchin} N.J. Hitchin: Twistor spaces, Einstein metrics and
isomonodromic deformations. J. Diff. Geom. \textbf{42} (1995) 30

\bibitem{Jeffrey}
L.C. Jeffrey: Extended moduli spaces of flat connections on Riemann
surfaces. Math. Ann. \textbf{298} (1994) 667

\bibitem{Kirillov}
A.A. Kirillov:
Elements of the theory of representations.
Springer Verlag 1976

\bibitem{Lerman-Sjamaar} E. Lerman, R. Sjamaar:  Stratified
symplectic spaces and reduction. Math. Ann. \textbf{134} (1991) 375

\bibitem{Lev-Soi}
S. Levendorskii, Y. Soibelman:
Algebras of functions on compact quantum groups,
Shubert cells and quantum tori.
Commun. Math. Phys. \textbf{139} (1991) 141

\bibitem{Lu} J.-H. Lu: Moment mappings and reduction of Poisson
actions. Proc. of the Sem. Sud-Rodanian de Geometrie \'{a} Berkeley
(1989) Springer Verlag MSRI Series 1991

\bibitem{Lu-Weinstein} J.-H. Lu, A. Weinstein: Poisson Lie groups,
dressing transformations and Bruhat decompositions. J. Diff. Geom.
\textbf{31}  (1990) 510

\bibitem{Semenov} M.A. Semenov-Tian-Shansky: What is a classical
$r$-matrix.  Funct. Anal. Appl. \textbf{17} (1983) 259

\bibitem{Sklyanin} E.K. Sklyanin:
Some algebraic structures related to the Yang-Baxter equation.
Funk. Anal. i ego prilozh. \textbf{16} (1982) 27

\bibitem{Weinstein} A. Weinstein:
The local structure of Poisson manifolds.
J.Diff.Geom. \textbf{18} (1983) 523

\end{thebibliography}
\end{document}